\documentclass[12pt]{article}
\usepackage{amssymb}
\usepackage{latexsym,bm}

\newtheorem{theorem}{Theorem}
\newtheorem{corollary}[theorem]{Corollary}

\newtheorem{lemma}[theorem]{Lemma}

\usepackage{amsmath}

\newcommand{\bea}{\begin{eqnarray*}}
\newcommand{\eea}{\end{eqnarray*}}
\newcommand{\be}{\begin{equation}}
\newcommand{\ee}{\end{equation}}
\newcommand{\ben}{\begin{eqnarray*}}
\newcommand{\een}{\end{eqnarray*}}

\newcommand{\C}{{\Bbb C}}

\date{}
\begin{document}
\title{Nonsymmetric normal entry patterns with the maximum number of distinct indeterminates}

\author{Zejun Huang\\
College of Mathematics and Econometrics,
 Hunan University\\
Changsha, Hunan 410082, China\\
{\tt mathzejun@gmail.com} \and Xingzhi Zhan\footnote{Corresponding author}\\
Department of Mathematics,
East China Normal University\\
Shanghai 200241, China\\
{\tt zhan@math.ecnu.edu.cn} }\maketitle

\begin{abstract}
We prove that a nonsymmetric normal entry pattern of order $n$ ($n\ge 3$)
has at most $n(n-3)/2+3$ distinct indeterminates and up to permutation similarity this number is attained by a unique pattern
which is explicitly described.
\end{abstract}

{\bf AMS classifications:} 05D99, 05B20, 15A99, 15B05

{\bf Keywords:} Entry pattern, normal matrix, 0-1 matrix

\section{Introduction}

Symmetric matrices, Toeplitz matrices, Hankel matrices and circulant matrices all require repetitions of some entries. These special matrices
suggest that we define a new concept for investigation of the general situation. Given a set $S,$ we denote by $M_n(S)$ the set of $n\times n$ matrices
whose entries are from $S.$ If $S=\{x_1,\ldots,x_k\}$ is a finite set, we write $M_n\{x_1,\ldots, x_k\}$ for $M_n(S).$

{\bf Definition 1} Let $x_1,x_2,\ldots,x_k$ be distinct indeterminates. We call a matrix in $M_n\{x_1,x_2,\ldots, x_k\}$ an {\it entry pattern.}

Thus an entry pattern is a matrix whose entries are indeterminates some of which may be equal. For example, among
$$
A=\left[\begin{array}{cc} x&y\\
                           y&z\end{array}\right],\,\,\,
B=\left[\begin{array}{cc} 2x&x+y\\
                           -z&w\end{array}\right],\,\,\,
C=\left[\begin{array}{cc} 3&x\\
                           y&z\end{array}\right]
$$
$A$ is an entry pattern while $B$ and $C$ are not. Rectangular entry patterns are defined similarly.

The spirit of entry patterns is that sometimes we can deduce properties of certain special matrices by just looking at the patterns
of their entries without knowing the actual entries. This is possible.
For example, every real symmetric matrix has all real eigenvalues and every complex circulant matrix is normal [4, p.5].
Entry patterns will serve the study of matrices over fields. To avoid unnecessary technical complications we consider only real matrices.
Given an entry pattern $A,$ we denote by $Q(A)$ the set of the real matrices obtained by specifying the values of the indeterminates of $A.$
Thus
$$
\left[\begin{array}{ccc}2&3&5\\
                        5&2&3\\
                        3&5&2\end{array}\right]\in Q(A)\,\,\,{\rm with}\,\,\,A=\left[\begin{array}{ccc}x&y&z\\
                                                                                                       z&x&y\\
                                                                                                       y&z&x\end{array}\right].
$$
Conversely,
$$
\left[\begin{array}{ccc}1&1&1\\
                        2&3&4\\
                        2&4&3\end{array}\right]\,\,\,{\rm and}\,\,\,
 \left[\begin{array}{ccc}5&5&5\\
                         6&7&8\\
                         6&8&7\end{array}\right]\,\,\,{\rm have}\,\,\,{\rm the}\,\,\,{\rm same}\,\,\,{\rm entry}\,\,\,{\rm pattern}
 \left[\begin{array}{ccc}x&x&x\\
                         y&z&w\\
                         y&w&z\end{array}\right].
$$

We denote by $A^T$ the transpose of a matrix $A.$ Recall that a real matrix $A$ is said to be {\it normal} if $AA^T=A^TA.$
Including symmetric matrices and orthogonal matrices as subclasses, normal matrices have nice properties and they are an important
topic in matrix analysis. See [1, Chapters VI and VII] and [2, Chapter 8].

{\bf Definition 2} A square entry pattern $A$ is said to be {\it normal} if every matrix in $Q(A)$ is normal.

Symmetric entry patterns are obviously normal. There are many nonsymmetric entry patterns. We will determine
the maximum number of distinct indeterminates in a nonsymmetric normal entry pattern of a given order and the patterns
that attain this number.

\section{Main results}

The main result is the following theorem.

\begin{theorem}
 Let $n\ge 3$ be an integer, and let $A$ be a nonsymmetric normal entry pattern of order $n$ with $k$ distinct entries. Then
 $k\le n(n-3)/2+3$, where equality holds if and only if $A$ is permutation similar to a pattern of the form
\begin{equation}\label{(1)}
\left[\begin{array}{cccc|ccc}
 x_{11}&x_{12}&\cdots&x_{1,n-3}&y_1&y_1&y_1\\
 x_{12}&x_{22}&\cdots&x_{2,n-3}&y_2&y_2&y_2\\
 \vdots&\vdots&\ddots&\vdots&\vdots&\vdots&\vdots\\
 x_{1,n-3}&x_{2,n-3}&\cdots&x_{n-3,n-3}&y_{n-3}&y_{n-3}&y_{n-3}\\
 \hline
 y_1&y_2&\cdots&y_{n-3}&z&u&v\\
 y_1&y_2&\cdots&y_{n-3}&v &z&u\\
  y_1&y_2&\cdots&y_{n-3}&u&v &z
  \end{array}\right]
 \end{equation}
 where $u,v,z$, $y_i$, $x_{ij}$, $1\le i\le j\le n-3$, are distinct indeterminates.
 \end{theorem}

Note that the matrix in (1) is of the form $\left[\begin{array}{cc}X&Y\\Y^T&Z\end{array}\right]$ where $X$ is symmetric, the entries in each
row of $Y$ are equal and $Z$ is a circulant matrix of order $3.$

It follows from Theorem 1 that the number of distinct indeterminates in a nonsymmetric normal entry pattern of order $n$ ($n\ge 3$) can be any number
in the interval $[2,\, n(n-3)/2+3].$

To prove Theorem 1 we need some lemmas. Denote by $J_n$ the $n\times n$ matrix with all entries equal to $1.$ An entry pattern $A$ in $M_n\{x_1,\ldots,x_k\}$
can be expressed uniquely as $A=\sum_{i=1}^kx_iA_i$ where $A_1,\ldots, A_k$ are 0-1 matrices with $\sum_{i=1}^kA_i=J_n.$ We call $A_i$ the {\it coefficient
matrix} of $x_i$ in $A.$ Here and in the sequel we view an entry pattern as a matrix whose entries are polynomials over the field of real numbers
${\mathbb R}$ so that addition and multiplication of entry patterns are defined in the usual way.

\begin{lemma}\label{lemma 2} Let $A_i$ be the coefficient matrix of $x_i$ in an entry pattern $A\in M_n\{x_1,\ldots,x_k\},$ $i=1,\ldots,k$. Then $A$ is a normal entry pattern if and only if each $A_i$ is normal and
\begin{equation}\label{(2)}
A_iA_j^T+A_jA_i^T=A_i^TA_j+A_j^TA_i\,\,\textrm{ for all }\,\,1\le i<j\le k.
\end{equation}
\end{lemma}

{\bf Proof.} Since $A=\sum_{i=1}^kx_iA_i$,  we have
$$AA^T=(\sum_{i=1}^kx_iA_i)(\sum_{j=1}^kx_jA_j)^T=\sum_{i=1}^kx_i^2A_iA_i^T+\sum_{1\le i<j\le k}x_ix_j(A_iA_j^T+A_jA_i^T)$$
and
$$
A^T A=(\sum_{i=1}^kx_iA_i)^T(\sum_{j=1}^kx_jA_j)=\sum_{i=1}^kx_i^2A_i^TA_i+\sum_{1\le i<j\le k}x_ix_j(A_i^T A_j +A_j^T A_i).
$$
Now suppose $A$ is normal. Then $AA^T=A^TA$ for any real values $x_1,\ldots,x_k.$ In this equality fixing any $i$ with $1\le i\le k$ and setting
$x_i=1$ and all other $x_j=0$ we obtain $A_iA_i^T=A_i^TA_i;$ i.e., $A_i$ is normal. Then fixing any pair $i,j$ with $1\le i<j\le k$ and setting
 $x_i=x_j=1$ and $x_t=0$ for all $t\in \{1,\ldots,k\}\setminus\{i,j\}$ we obtain (2).

 The converse implication is obvious. $\Box$

 Lemma 2 and its proof show that an entry pattern $A$ in $M_{n}\{x_1,\ldots,x_k\}$ is normal if and only if $A$ is normal for all $x_1,\ldots,x_k\in \{0,1\}$.

Let $B=(b_{ij})$ be a normal 0-1 matrix of order $n.$  For $1\le i\le n$, denote by $r_i$ and $c_i$ the $i$-th row sum and the $i$-th column sum of $B$ respectively, and denote by  $r'_i$  ($c'_i$) the sum of off-diagonal entries in the $i$-th row  (column) of $B$. Then $r'_i=r_i-b_{ii},$ $c'_i=c_i-b_{ii}$
for $1\le i\le n$. Equating the $i$-th diagonal entries of both sides of $BB^T=B^TB$ we have
\begin{equation}\label{(3)}
r_i=c_i\quad{\rm and}\quad r'_i=c'_i,\quad i=1,\ldots,n.
\end{equation}

\begin{corollary}\label{Corollary 3}
Let $A_1$ be the coefficient matrix of $x_1$ in $A\in M_n\{x_1,x_2\}.$ Then $A$  is a normal entry pattern if and only if $A_1$ is normal.
\end{corollary}

{\bf Proof.} If $A$ is a normal entry pattern, then by Lemma 2 $A_1$ is normal. Conversely suppose $A_1$ is normal. The equalities $r_i=c_i,$
$i=1,\ldots,n$ in (3) imply $J_nA_1^T+A_1J_n=J_nA_1+A_1^TJ_n,$ from which it follows that $A_2$ is normal and
$A_1A_2^T+A_2A_1^T=A_1^TA_2+A_2^TA_1$ where we have used the fact that $A_2=J_n-A_1.$ Applying Lemma 2 again in another direction we conclude that $A$ is normal. $\Box$

We remark that there is no known characterization of normal 0-1 matrices; see [3].

Throughout we denote by $f(B)$ the number of ones in a 0-1 matrix $B,$ by $O_k$ the zero matrix of order $k$ and by $I_k$ the identity matrix of order $k.$
Sometimes we omit the subscript $k$ if the order is clear from the context.
For square matrices $A,B$ the notation $A\oplus B$ means the block diagonal matrix ${\rm diag}(A,B).$ The notation $\equiv$ means that we denote something.

\begin{lemma}\label{Lemma 4}
Let $B$ be an $n\times n$ normal 0-1 matrix with $n\ge 2.$ Then
 \begin{itemize}
 \item[(i)] $f(B)=1$ if and only if $B$ is permutation similar to $1\oplus O_{n-1}$.
 \item[(ii)] $f(B)=2$ if and only if $B$ is permutation similar to $I_2\oplus O_{n-2}$ or $\begin{bmatrix}0&1\\1&0\end{bmatrix}\oplus O_{n-2}$.
\item[(iii)] $f(B)=3$ if and only if $B$ is permutation similar to one of the following four matrices:
\begin{equation}\label{(4)}
(a)~~I_3\oplus O_{n-3},~~(b)~~\begin{bmatrix}1&1\\1&0\end{bmatrix}\oplus O_{n-2},~~(c)~~ \begin{bmatrix}1&0&0\\0& 0&1\\0&1&0\end{bmatrix}\oplus O_{n-3},~~ (d)~~\begin{bmatrix}0&1&0\\0& 0&1\\1&0&0\end{bmatrix}\oplus O_{n-3}.
\end{equation}
 \end{itemize}
\end{lemma}

{\bf Proof.} The sufficiency is clear. Next we prove the necessity.

If $f(B)=1$, by (3) the only nonzero entry must be a diagonal entry. Hence $B$ is permutation similar to  $1\oplus O_{n-1}.$

If $f(B)=2$ and one diagonal entry is nonzero, then by (3) the other nonzero entry  of $B$ is also a diagonal entry. Thus $B$ is permutation similar to $I_2\oplus O_{n-2}$.

 If $f(B)=2$ and the diagonal entries are all zero, then by (3) $B$ is permutation similar to
 $\begin{bmatrix}0&1\\1&0\end{bmatrix}\oplus O_{n-2}$.

If $f(B)=3$, we distinguish three cases according to the number of nonzero diagonal entries.
By (3) we know that $B$ cannot have exactly two nonzero diagonal entries. If $B$ has three nonzero diagonal entries, then all the off-diagonal entries of $B$ are zero and $B$ is permutation similar to the matrix in $(a)$. If $B$ has only one nonzero diagonal entry, then $B$ must be symmetric. Hence it is permutation similar to $(b)$ or $(c)$. If $B\equiv(b_{ij})$ has no nonzero diagonal entry, by (3) we know that each row (column) of $B$ has at most one nonzero entry. Using a permutation similarity transformation if necessary we may assume that the first row of $B$ has an entry equal to one, say, $b_{12}=1$. Then there is an entry $b_{j1}=1$. If $j=2$, then $f(B)=3$ forces the third nonzero entry in $B$ to be a diagonal entry,  a contradiction.  Hence $j\not\in \{1,2\}$ and we may assume $j=3$, since we can permute rows $j$ and $3$ and then permute the columns $j$ and $3$ if necessary. Now we have proved that $B$ is permutation similar to
$$C=\begin{bmatrix}0&1&0&0&\cdots&0\\
                 0&0&?&?&\cdots&?\\
                 1&0&0&0&\cdots&0\\
                 0&0&?&?&\cdots&?\\
                 \vdots&\vdots&\vdots&\vdots&\cdots&\vdots\\
                 0&0&?&?&\cdots&?
                 \end{bmatrix}\equiv(c_{ij}).$$
Again, by (3) we know that the third nonzero entry in $C$ is in the second row and in the third column, which means $c_{23}=1$. Hence $B$ is permutation similar to $(d)$. $\Box$

\begin{lemma}\label{Lemma 5}
Let  $B=\begin{bmatrix}B_{1}&B_{2}\\B^T_{2}&B_{3}\end{bmatrix}$ be a square real matrix or an entry pattern. If $B_{1}$ is symmetric, 
then $B$ is normal if and only if  $B_{3}$ is normal and $B_{2}B_{3}=B_{2}B_{3}^T$.
\end{lemma}

{\bf Proof.} Since
$$BB^T=\begin{bmatrix}B^2_{1}+B_{2}B^T_{2}&B_{1}B_{2}+B_{2}B_{3}^T\\
B_{2}^TB_{1}+B_{3}B_{2}^T&B_{2}^TB_{2}+B_{3}B_{3}^T\end{bmatrix}$$
and
$$B^T B=\begin{bmatrix}B^2_{1}+B_{2}B^T_{2}&B_{1}B_{2}+B_{2}B_{3}\\
B_{2}^TB_{1}+B^T_{3}B_{2}^T&B_{2}^TB_{2}+B_{3}^T B_{3}\end{bmatrix},$$
the conclusion is clear. $\Box$

\begin{lemma}\label{Lemma 6}
Let $A$ be a normal entry pattern in $M_n\{x_1,\ldots,x_k\}$ and let $A_i$ be the coefficient matrix of $x_i$ in $A.$
\begin{itemize}
\item[(i)] If $A_i$ has exactly one nonzero entry, then $A$ is permutation similar to
$$
\begin{bmatrix}x_i&a\\a^T&B\end{bmatrix}.
$$
\item[(ii)] If $A_i$ has exactly two nonzero entries, then $A$ is permutation similar to
\begin{equation}\label{(5)}
\begin{bmatrix}x_i&x_j&b\\x_j&x_i&c\\b^T&c^T&B\end{bmatrix}\quad\textrm{ or }\quad \begin{bmatrix}x_j&x_i&b\\x_i&x_k&c\\b^T&c^T&B\end{bmatrix}
\end{equation}
where $j\ne i$ and $k\ne i$.
\end{itemize}
\end{lemma}

{\bf Proof.}
(i)  By Lemma 2, $A_i$ is a normal 0-1 matrix with only one nonzero entry. By Lemma 4, $A_i$ is permutation similar to $1\oplus O_{n-1}$.  Without loss of generality we assume $A_i=1\oplus O_{n-1}$. For $j\in \{1,\ldots,k\}\setminus \{i\}$, partition $A_j$ as
$$
A_j=\begin{bmatrix}0&b_j\\
                  c_j^T&B_j \end{bmatrix}.
$$
Then
 $$A_iA_j^T+A_jA_i^T=\begin{bmatrix}0&c_j\\c_j^T&O\end{bmatrix} \quad
\textrm{ and }\quad
 A_i^T A_j+A_j^T A_i=\begin{bmatrix}0&b_j\\b_j^T&O\end{bmatrix}.$$
By (2) we have $b_j=c_j$. Hence $A=\sum_jx_jA_j$ has the required form in (i).

(ii) Applying Lemma 4, we may assume $A_i=I_2\oplus O_{n-2}$ or
$\begin{bmatrix}0&1\\1&0\end{bmatrix}\oplus O_{n-2}$. If $A_i=I_2\oplus O_{n-2}$, for $j\in \{1,\ldots,k\}\setminus \{i\}$, partition $A_j$ as
$$A_j=\begin{bmatrix}0&\lambda_j&b_j\\
                      \theta_j&0&c_j\\
                      r_j^T&s_j^T&B_j        \end{bmatrix}$$
                      where $\lambda_j,\theta_j\in\{0,1\}.$
Using the same arguments as above, we have $r_j=b_j,s_j=c_j$ for all $j\ne i$.  Moreover, by (3) we have $\lambda_j=\theta_j$ for all $j\ne i$.
Therefore, $A$ is of the first form in (5).

If $A_i=\begin{bmatrix}0&1\\1&0\end{bmatrix}\oplus O_{n-2}\equiv P\oplus O_{n-2}$, then for $j\in \{1,\ldots,k\}\setminus \{i\}$, we partition $A_j$ as
 $$A_j=\begin{bmatrix}
B_{j1}&B_{j2}\\B_{j3}&B_{j4}
\end{bmatrix}$$
with $B_{j1}$ a $2\times 2$ matrix. Since $\sum_tA_t=J_n$, $B_{j1}$ is diagonal. We have
$$A_iA_j^T+A_jA_i^T=\begin{bmatrix}PB_{j1}+B_{j1}P&PB_{j3}^T\\B_{j3}P&O\end{bmatrix}$$
and
$$A_i^T A_j+A_j^T A_i=\begin{bmatrix}PB_{j1}+B_{j1}P&PB_{j2}\\B_{j2}^T P&O\end{bmatrix}.$$
By (2) we have $B_{j3}=B_{j2}^T$. Hence  $A=\sum_tx_tA_t$ has the second  form in (5). $\Box$

\begin{lemma}\label{Lemma 7}
 A  $2\times 2$  entry pattern is normal if and only if it is symmetric.
\end{lemma}

{\bf Proof.}
It is clear that any symmetric entry pattern is normal.
If a $2\times 2$ entry pattern $A$ is normal,  Lemma 2 implies that each coefficient matrix $A_i$ is a $2\times 2$ normal 0-1 matrix.
Then (3) implies that each $A_i$ is symmetric and hence $A$ is symmetric. $\Box$

\begin{lemma}\label{Lemma 8}
Theorem 1 is true for $n=3$.
\end{lemma}

{\bf Proof.}
By a direct computation one can verify that any entry pattern of the form
 \begin{equation}\begin{bmatrix}z&u&v\\ v&z&u\\u&v&z\end{bmatrix}\end{equation}
 is normal.

 Conversely, let $A$ be a nonsymmetric normal entry pattern of order 3.
 Suppose $A$ has at least 4 distinct entries. Since $A$ has only 9 entries, there is an entry $x_i$, say, $x_1$, which
  appears exactly once or twice in $A$, i.e., $f(A_1)\le 2$.
 If $f(A_1)=1$, then by  Lemma 6(i), $A$ is permutation similar to
 $$\begin{bmatrix}x_1&a\\a^T&B\end{bmatrix}.$$
 Moreover, by Lemma 5, $B$ is a $2\times 2$ normal entry pattern, which is symmetric by Lemma 7. Hence $A$ is symmetric, a contradiction.
 If $f(A_1)=2$, by  Lemma 6(ii), $A$ is symmetric, a contradiction.  Hence $A$ has at most 3 distinct entries. Note that with $n=3,$
 $n(n-3)/2+3=3.$

 Suppose $A$ has exactly 3 distinct entries $x_1,x_2,x_3$. If one of $x_1,x_2,x_3$ appears exactly once or twice in $A$,
 then using the same arguments as above we deduce that $A$ is symmetric. Hence we have $f(A_1)=f(A_2)=f(A_3)=3$.
 By Lemma 2, $A_1,A_2,A_3$ are all normal 0-1 matrices. Applying Lemma 4(iii) we conclude that each of $A_1,A_2$ and $A_3$ is permutation
 similar to one of the following four matrices:
 $$
 I_3,\quad\begin{bmatrix}1&1&0\\1&0&0\\0&0&0\end{bmatrix},\quad\begin{bmatrix}1&0&0\\0&0&1\\0&1&0\end{bmatrix},\quad
  P=\begin{bmatrix}0&1&0\\0& 0&1\\1&0&0\end{bmatrix}.
 $$
 Note that the first three matrices are symmetric. Since $A$ is nonsymmetric, one of $A_1,A_2$ and $A_3$ is permutation similar to $P.$
 Since $A_1+A_2+A_3=J_3,$ by considering the diagonal entries we deduce that one of the other two coefficient matrices is $I_3.$
 Thus, there is a permutation matrix $Q$ such that
$$\{A_1,\,A_2,\,A_3\}=\{I_3,\,Q^TPQ,\,Q^TP^TQ\}.$$
It follows that $A$ is permutation similar to a pattern of the form (6). $\Box$

\par
{\bf  Proof of Theorem 1.}
 Denote by $\phi(G)$ the number of distinct entries in an entry pattern $G$. First we use induction on the order $n$ to prove that if $A$ is a nonsymmetric normal entry pattern of order $n$, then $\phi(A)\le n(n-3)/2+3$. Lemma 8 shows that this is true for $n=3.$ Now let $n\ge 4$ and assume that
 the conclusion is true for all entry patterns of order $n-1$. Let $A$ be a nonsymmetric normal entry pattern of order $n.$

 To the contrary, assume that $\phi(A)\ge n(n-3)/2+4$. Then $A$ contains an entry, say, $x_1$, which appears exactly once or twice in $A$. Otherwise, each entry appears at least 3 times and we have
$$3[n(n-3)/2+4]>n^2,$$
which contradicts the fact that $A$ has only $n^2$ entries. Let $A_i$ be the coefficient matrix of $x_i$ in $A.$

If $f(A_1)=1$, then by Lemma 6(i) $A$ is permutation similar to
\begin{equation}\label{(7)}
A^{(1)}= \begin{bmatrix}x_1&a\\a^T&B\end{bmatrix}.
\end{equation}
Since $A$ is nonsymmetric, so is $B.$ By Lemma 5, $B$ is a normal entry pattern of order $n-1.$ Applying the induction hypothesis to $B,$
we have $\phi(B)\le (n-1)(n-4)/2+3$. Denote by $\theta(a)$ the number of those distinct entries in $A^{(1)}$ that appear only in $a$ and $a^T.$
Then
\begin{eqnarray*}
\theta(a)&=&\phi(A)-\phi(B)-1\\
         &\ge&[n(n-3)/2+4]-[(n-1)(n-4)/2+3]-1\\
         &=&n-2.
\end{eqnarray*}
Since $a$ has $n-1$ components, it follows that there are at least $n-3$ distinct entries in $a$ that appear exactly twice in $A^{(1)}$.
Without loss of generality, we assume $f(A_2)=\cdots=f(A_{n-2})=2.$ Lemma 6(ii) indicates that for each $2\le i\le n-2,$ the two rows of $A^{(1)}$
in which the two $x_i's$ lie and the corresponding two columns are symmetric. Thus, $A^{(1)}$ is permutation similar to

\begin{equation}\label{(8)}
A^{(2)}=\left[ \begin{array}{c|ccc|c}
x_1&x_2&\cdots&x_{n-2}&b\cr \hline
x_2& & & &\cr
\vdots&&C&&E\cr
x_{n-2}&&&&\cr \hline
b^T&&E^T&&F
 \end{array}\right]
\end{equation}
where $C$ is symmetric. By Lemma 5, $F$ is a $2\times 2$ normal entry pattern which must be symmetric by Lemma 7. Hence $A^{(2)}$ is symmetric
and consequently $A$ is symmetric, a contradiction.

 If $f(A_1)=2$, then by Lemma 6(ii), $A$ is permutation similar to a matrix of one of the two forms in (5)
 with $i=1$. Repartition
 \begin{equation}\label{(9)}
  \begin{bmatrix}x_1&x_j&b\\x_j&x_1&c\\b^T&c^T&B\end{bmatrix}
 =\begin{bmatrix}x_1&a\\a^T&S\end{bmatrix}, \quad \begin{bmatrix}x_j&x_1&b\\x_1&x_k&c\\b^T&c^T&B\end{bmatrix}
 =\begin{bmatrix}x_j&a\\a^T&S\end{bmatrix}.
 \end{equation}
In both cases we have $\theta(a)\ge n-2$. Using the same arguments as above we deduce that $A$ is symmetric, a contradiction.

Therefore, $\phi(A)\le n(n-3)/2+3$.

Next we use induction on the order $n$ to prove that if $A$ is a nonsymmetric normal entry pattern of order $n$ with  $\phi(A)=n(n-3)/2+3$, then $A$ is permutation similar to a pattern of the form (1). Lemma 8 shows that this is true for $n=3.$ Now let $n\ge 4$ and assume that
the conclusion is true for all entry patterns of order $n-1$. Let $A$ be a nonsymmetric normal entry pattern of order $n$ with $\phi(A)=n(n-3)/2+3.$

From $4[n(n-3)/2+3]> n^2$ we know that there is at least one entry that appears less than 4 times in $A$. Suppose $x_1$ is an entry that appears the least times in $A$. Then $f(A_1)\le 3$. We distinguish three cases.

{\it Case 1.} $f(A_1)=1$. By Lemma 6(i) $A$ is permutation similar to a pattern of the  form (7).
Applying Lemma 5 we deduce that the matrix $B$ in (7) is a normal entry pattern of order $n-1.$ By what we have proved above,
$$\phi(B)\le(n-1)(n-4)/2+3.$$
If $\phi(B)<(n-1)(n-4)/2+3$, then  there are at least
$$\phi(A)-\phi(B)-1\ge n-2$$
distinct entries in $a$ that do not appear in $B$. Using the same arguments as above, we conclude that $A$ is symmetric, a contradiction. Hence we have  $\phi(B)=(n-1)(n-4)/2+3$. By the induction hypothesis, $B$ is permutation similar to a matrix of the  form (1) and hence $A$ is permutation similar to
$$
H=\begin{bmatrix}
 x_{1}&c_{2}&\cdots&c_{n-3}&c_{n-2}&c_{n-1}&c_{n}\\
c_2&x_{22}&\cdots&x_{2,n-3}&y_2&y_2&y_2\\
 \vdots&\vdots&\ddots&\vdots&\vdots&\vdots&\vdots\\
c_{n-3}&x_{2,n-3}&\cdots&x_{n-3,n-3}&y_{n-3}&y_{n-3}&y_{n-3}\\
 c_{n-2}&y_2&\cdots&y_{n-3}&z&u&v\\
c_{n-1}&y_2&\cdots&y_{n-3}&v &z&u\\
  c_{n}&y_2&\cdots&y_{n-3}&u&v &z
 \end{bmatrix}.
$$
Suppose there are exactly $d$ distinct entries among $c_2,\ldots,c_n$ that do not appear in $B$. Then
$$1+d+\phi(B)=n(n-3)/2+3,$$
 yielding $d=n-3$. Now it suffices to prove that $c_{n-2}=c_{n-1}=c_{n}$, which forces $c_2,\ldots,c_{n-3}$ to be distinct entries, and hence $H$ has the required form.

Let $y_i=x_{ij}=0$ for all $2\le i\le j\le n-3$. Then $H$ is permutation similar to
\begin{equation*}
 K=\begin{bmatrix}
 & & &c_{2}& &&\\
  &  & &\vdots && & \\
 & & &c_{n-3}& & & \\
c_{2}&\cdots&c_{n-3}&x_1&c_{n-2}&c_{n-1}&c_{n}\\
&&&c_{n-2} &z&u&v\\
&&&c_{n-1} &v &z&u\\
 &&& c_{n} &u&v &z
 \end{bmatrix}\equiv\begin{bmatrix}O&p^T&O\\p&x_1&q\\O&q^T&F\end{bmatrix}
 \end{equation*}
 where $p=(c_2,\ldots,c_{n-3}),$ $q=(c_{n-2},c_{n-1},c_{n})$.
 Now $KK^T=K^TK$ implies $qF=qF^T$. Setting $z=v=0$  and $u=1$ in $F$ we get
 $$0=q(F-F^T)=(c_n-c_{n-1},\,c_{n-2}-c_n,\,*),$$
which gives $c_{n-2}=c_{n-1}=c_{n}$.

{\it Case 2.}  $f(A_1)=2$. We will show that this case cannot happen. By Lemma 6(ii), $A$ is permutation similar to one of the two forms in  (5). In both cases, we can repartition $A$ as (9). Applying Lemma 5 we deduce that the matrix $S$ in (9) is a nonsymmetric normal entry pattern of order $n-1.$ 
By what we have proved above,
$$\phi(S)\le(n-1)(n-4)/2+3.$$
If $\phi(S)<(n-1)(n-4)/2+3$, then
$$\theta(a)\ge\phi(A)-\phi(S)-1\ge n-2.$$ Using the same arguments as above, we deduce that $A$
 is permutation similar to a matrix of form (8), which is symmetric, a contradiction. Hence we have  $\phi(S)=(n-1)(n-4)/2+3$. 
 By the induction hypothesis, $S$ is permutation similar to a matrix of the  form (1) and hence $A$ is permutation similar to
  $$
\begin{bmatrix}
 w&c_{2}&\cdots&c_{n-3}&c_{n-2}&c_{n-1}&c_{n}\\
c_2&x_{22}&\cdots&x_{2,n-3}&y_2&y_2&y_2\\
 \vdots&\vdots&\ddots&\vdots&\vdots&\vdots&\vdots\\
c_{n-3}&x_{2,n-3}&\cdots&x_{n-3,n-3}&y_{n-3}&y_{n-3}&y_{n-3}\\
 c_{n-2}&y_2&\cdots&y_{n-3}&z&u&v\\
c_{n-1}&y_2&\cdots&y_{n-3}&v &z&u\\
  c_{n}&y_2&\cdots&y_{n-3}&u&v &z
 \end{bmatrix}
 $$
 where $w=x_1$ or $w=x_j.$ As in Case 1, we can prove that $c_{n-2}=c_{n-1}=c_n.$  Since this matrix has $n(n-3)/2+3$ distinct entries by assumption ,
 it follows that each of the diagonal entries $w,x_{22},\ldots,x_{n-3,n-3}$ appears exactly once, which contradicts 
 the assumption that $x_1$ appears the least time $2.$

 {\it Case 3.} $f(A_1)=3$. We will show that this case cannot happen. If $n\ge 7$, then $3[n(n-3)/2+3]>n^2$, a contradiction. Hence we have  $n\le 6$. Moreover, if $n=6$, each of the $n(n-3)/2+3$ distinct entries appears exactly 3 times. If $n=4$ or $n=5$, then one of the $n(n-3)/2+3$ distinct entries  appears exactly 4 times and the others appear exactly 3 times. By Lemma 4(iii), we may assume $A_1$ has one of the forms in (4).

{\it Subcase 3.1.} $A_1$ has form (b) in (4). For any $2\le j\le n(n-3)/2+3$,  partition $A_j$ as
 \begin{equation}\label{(10)}
 A_j=\begin{bmatrix}B_{j1}&B_{j2}\\B_{j3}&B_{j4}\end{bmatrix}
 \end{equation}
with $B_{j1}$ a $2\times 2$ matrix. Applying (2) we deduce that $B_{j3}=B_{j2}^T$ and hence $A$ has the form
\begin{equation}\label{eq11}
A=\begin{bmatrix}C_1&C_2\\C_2^T&C_3 \end{bmatrix}
\end{equation}
with $C_1=\begin{bmatrix}x_1&x_1\\x_1&\beta\end{bmatrix}$. By Lemma 5 and the assumption on $A,$  $C_3$ is a nonsymmetric normal entry pattern.

If $n=4$, the $2\times 2$ normal pattern $C_3$ must be symmetric, a contradiction.

If $n=5$, $C_3$  contains at least 3 distinct entries. By Lemma 8, $C_3$ is permutation similar to
\begin{equation}\label{(12)}
\begin{bmatrix}z&u&v\\v&z&u\\u&v&z\end{bmatrix}.
\end{equation}
Note that since $C_2$ and $C_2^T$ are in symmetric positions in $A,$ the only way for an entry in $C_2$ to appear exactly 3 times in $A$ is
appearing once in a diagonal position. Since all diagonal entries in $C_3$ are equal, then at most two entries in $C_2$ appear
 exactly 3 times in $A$. Thus at least two entries in $C_2$ appear an even number of times, a contradiction.

If $n=6$, it is impossible for all the 8 entries in $C_2$ to appear exactly 3 times in $A$, since $A$ has only 6 diagonal positions.

{\it Subcase 3.2.} $A_1$ has form (a), (c)  or (d) in (4). Again,   for any $2\le j\le n(n-3)/2+3$,  partition $A_j$ as form (10) with $B_{j1}$ a $3\times 3$ matrix. Applying (2) we deduce that $A$ has   form (11) with $C_1$ a $3\times 3$ entry pattern.

If $n=4$ or $n=5$, then by (3), $C_3$ is symmetric, and hence $C_1$ is a nonsymmetric normal entry pattern by Lemma 5.
$C_1$ has 3 distinct entries.  Applying Lemma 8, we deduce that $C_1$ is permutation similar to a pattern of form (12).
If $n=4$, then $C_2$ has at least one entry, say, $x_2$, appears exactly 3 times in $A$. It follows that $C_3=x_2$ and $A_2$ is permutation similar to $\begin{bmatrix}1&1\\1&0\end{bmatrix}\oplus O_2$. Applying subcase 3.1 we get a contradiction. If $n=5$, then $C_2$ has 6 entries
 at least 4 of which appear exactly 3 times in $A$, which is impossible since $C_1$ has three identical diagonal entries and $C_3$ has only 2 diagonal entries.

If $n=6$, then all the 9 entries in $C_2$  appear exactly 3 times in $A$, which is impossible since $A$ has only 6 diagonal positions.

So far we have proved that a nonsymmetric normal entry pattern of order $n$ with $n(n-3)/2+3$ distinct entries is permutation similar to a pattern
of the form (1). It remains to verify that an entry pattern $A$ of the form (1) is normal. Partition $A$ as $\left[\begin{array}{cc}X&Y\\Y^T&Z\end{array}\right]$ where $X$ is symmetric, the entries in each
row of $Y$ are equal and $Z$ is a circulant matrix of order $3.$
A direct computation shows that $AA^T=A^TA$.  This completes the proof. $\Box$

\section*{Acknowledgement}

The research of Huang was supported by the NSFC grant 11401197 and a Fundamental Research
Fund  for the Central Universities. This work was done when
Huang was visiting East China Normal University in Spring 2015. He thanks ECNU for its hospitality and support.
Zhan's research was supported by the Shanghai SF grant 15ZR1411500 and the NSFC grant 11371145.

\end{document}